\input amstex 
\documentstyle{amsppt}
\magnification1200
\NoRunningHeads
\refstyle{A}


\let\la=\langle 
\let\ra=\rangle 

\define\C{{\Bbb C}}
\define\Z{{\Bbb Z}}
\define\N{{\Bbb N}}
\def\part{\operatorname{part}}
\define\wt{\operatorname{wt}}

\define\lspan{\operatorname{\Bbb C\text{-span}}}
\define\ch{\operatorname{ch}}


\topmatter
\title{Some crystal Rogers-Ramanujan type identities}
\endtitle
\author Mirko Primc
\endauthor
\address
Univ\. of Zagreb,
Dept\. of Mathematics,
Bijeni\v{c}ka 30, Zagreb,
Croatia
\endaddress
\email
primc\@cromath.math.hr
\endemail
\thanks Partially supported by the Ministry of Science 
and Technology of the Republic of Croatia, grant 037002.
  \endthanks
\keywords
Rogers-Ramanujan identities, colored
partitions, partition ideals, perfect crystals,
affine Lie algebras, basic modules 
  \endkeywords
\subjclass Primary 05A19;
Secondary 17B37, 17B67
\endsubjclass
\abstract
  By using the KMN${}^2$ crystal base character formula  
for the basic $A_2^{(1)}$-module,
and the principally specialized Weyl-Kac character formula,
 we obtain a Rogers-Ramanujan type 
combinatorial identity for colored partitions. The difference 
conditions between parts are given by
the energy function of certain perfect $A_2^{(1)}$-crystal.
We also recall some other identities for this type of colored partitions,
but coming from the vertex operator constructions 
and with no apparent connection to the crystal base theory. 
\endabstract

\endtopmatter


\document

\heading 1. Introduction
\endheading

J. Lepowsky and R. Wilson gave in \cite{LW} a Lie-theoretic interpretation of 
Rogers-Rama\-nu\-jan identities in terms of representations of 
affine Lie algebra $\tilde\goth g=\goth {sl}(2,\Bbb C)\sptilde$.
The product sides of Rogers-Ramanujan identities follow from the
principally specialized Weyl-Kac character formula
for level 3 standard $\tilde\goth g$-modules, the sum sides follow from the
vertex operator construction of bases of level 3 standard $\tilde\goth g$-modules,
parameterized by partitions satisfying difference 2 conditions.

The Lepowsky-Wilson approach is also possible for other affine Lie algebras
and for other constructions of vertex operators, 
various combinatorial consequences are illustrated by
constructions given, for example, in \cite{C3}, \cite{LP}, \cite{Ma} and  \cite{Mi}.

In \cite{P1} there is a construction of the basic $A^{(1)}_\ell$-module
based on the Frenkel-Kac vertex operator formula. 
It was noted there that the combinatorial difference conditions 
arising from the vertex operator formula
coincide with the energy function in the
construction by paths of the basic representation of $A^{(1)}_{\ell -1}$ 
given in \cite{DJKMO}. In \cite{P2} this combinatorial connection for
basic modules is extended to other  classical affine Lie algebras, this time by
using more general crystal base character 
formula for standard modules due to S.-J. Kang, M. Kashiwara,
K. C. Misra, T. Miwa, T. Nakashima and A. Nakayashiki \cite{KMN${}^2$}. 

This combinatorial connection of the crystal base theory and the vertex operator
constructions suggested that it might be interesting to study Rogers-Ramanujan type
combinatorial identities for colored partitions, where differences are given by
energy functions of perfect crystals. So we start, roughly speaking, 
with colored partitions
$(n_1)_{\beta_1}\geqslant (n_2)_{\beta_2}\geqslant 
\dots \geqslant (n_s)_{\beta_s}>0$, where each 
number $n_r$ is ``colored'' with a ``color'' $\beta_r$ from the set of
nine ``colors'' $\{1,\dots, 9\}$. Analogous to the Rogers-Ramanujan case, 
we consider colored partitions satisfying difference conditions
$$
(n_{r})_{\beta_{r}}\geqslant (n_{r+1})_{\beta_{r+1}}+E_{\beta_r\beta_{r+1}} \ ,
$$
where differences $E_{\beta_r\beta_{r+1}}\in\{0,1,2\}$ are the values of an
energy function of certain perfect $\goth {sl}(3,\Bbb C)\sptilde$-crystal.
We obtain an identity for such colored partitions (Theorem 2.1) by
using the principally specialized Weyl-Kac character formula and the crystal 
base character formula \cite{KMN${}^2$}.

In the last section we recall some other identities for this type of colored partitions,
but coming from the vertex operator construction \cite{MP2} and with no apparent
connection to the crystal base theory.

I am grateful to Jim Lepowsky and Arne Meurman for many ideas and results 
implicit in this work, and to Ivica Siladi\' c for testing
the identities numerically.
I thank the ICTP in Trieste, where this work
started, for their kind hospitality.


\heading 2. A combinatorial identity
\endheading

Let $A$ be a nonempty set and denote by $\Cal P (A)$ the set
of all maps
$\pi : A \rightarrow \Bbb N$, where $\pi(a)$ equals zero for
all but finitely
many $a \in A$. Clearly $\pi$ is determined by its values
$(\pi(a) \mid a \in
A)$ and we shall also write $\pi$ as a monomial
$$
\pi = \prod_{a\in A} a^{\pi(a)}. 
$$
We shall say that $\pi$ is a partition and for $\pi(a)
> 0$ we shall say that
$a$ is a part of $\pi$. We define the length
$\ell(\pi)$ of $\pi$ by
$
\ell(\pi) = \sum_{a\in A} \pi(a).
$
We consider elements of
$A$ as partitions
of length 1, i.e. $A \subset \Cal P (A)$.
For $\rho,\pi \in \Cal P (A)$ we write 
 $\pi \supset \rho$ if $\pi (a)\geqslant\rho (a)  $ for all $a\in A$ and
 we say that  $\pi$ contains $\rho$.

 For 
$\rho,\pi \in \Cal P (A)$
we define $\pi \rho$ in $\Cal P (A)$ by
$(\pi\rho)(a)=\pi(a)+\rho(a)$,
$a \in A$. 
We shall say that $1=\prod_{a\in A} a^0$ is the partition
with no parts
 and length 0. Clearly 
$\Cal P (A)$ is a monoid. 
For lack of a better terminology, we shall say that $\Cal I
\subset \Cal P (A)$
is an ideal in the monoid $\Cal P (A)$ if $\rho \in
\Cal I$ and $\pi \in \Cal P (A)$
implies $\rho \pi \in \Cal I$. For such an $\Cal I$ we call the difference of sets
$\Cal P (A)\backslash \Cal I$
a partition ideal in $\Cal P (A)$ (cf. \cite{A1, Chapter 8}). Later on
we shall consider an ideal $\Cal I=\Cal P(A)\Cal D$ generated by a set $\Cal D$,
and the corresponding partition ideal we shall denote as
$$
\Cal P_{\Cal D}=\Cal P (A)\backslash (\Cal P(A)\Cal D).
$$
Let $a\mapsto |a|$ be a map from $A$ to $\N$.
Then we define the degree $\vert
\pi \vert$ of $\pi$ by
$$
\vert \pi \vert = \sum_{a\in A} |a|\,\pi(a),
$$
and we say that a part $a$ of $\pi$ has the degree $|a|$.

Now let 
$$
\Gamma =\{ 1,2,3,4,5,6,7,8,9\},
$$
we shall think of $\Gamma$ as a set of colors, and let $A=\Gamma\times\Z_{<0}$.
We shall write $(\alpha, i)=i_\alpha$,
$$
A=\{i_\alpha\mid  \alpha\in\Gamma, i\in\Z_{<0}\},
$$
and we shall think that $i_\alpha$ has a color $\alpha$.
We define a map $(-i)_\alpha\mapsto |(-i)_\alpha|$ for $i>0$ by
$$\alignedat3
&(-i)_1\mapsto 3i-2, \qquad\qquad &&(-i)_4\mapsto 3i, 
\qquad\qquad &&(-i)_9\mapsto 3i+2,\\
&(-i)_2\mapsto 3i-1, \qquad\qquad &&(-i)_5\mapsto 3i, 
\qquad\qquad &&(-i)_8\mapsto 3i+1,\\
&(-i)_3\mapsto 3i-1, \qquad\qquad &&(-i)_6\mapsto 3i, 
\qquad\qquad &&(-i)_7\mapsto 3i+1.
\endalignedat\tag{2.1}$$
So, for example, $(-5)_1$ has the color $1$ and the degree $|(-5)_1|=13$.
In general, we can think of $\pi\in\Cal P(A)$ as a colored partition of
the nonnegative integer $n=|\pi|$.

Let $E=(E_{\alpha\beta})_{\alpha,\beta\in\Gamma}$ be a matrix
$$
E =\left(\matrix 
2 & 2& 2& 1& 2& 2& 2& 2& 2\\
1 & 2& 1& 1& 2& 1& 2& 2& 2\\
1 & 1& 2& 1& 1& 2& 2& 2& 2\\
1 & 1& 1& 0& 1& 1& 1& 1& 1\\
0 & 0& 1& 1& 0& 1& 1& 2& 2\\
0 & 1& 0& 1& 1& 0& 2& 1& 2\\
0 & 1& 0& 1& 1& 0& 2& 1& 2\\
0 & 0& 1& 1& 0& 1& 1& 2& 2\\
0 & 0& 0& 1& 0& 0& 1& 1& 2 
\endmatrix\right) .
$$
Define a subset $\Cal D$ in $\Cal P(A)$ by
$$
\Cal D =\{i_\alpha i_\beta\mid  E_{\alpha\beta} E_{\beta\alpha}\geqslant 1\}
\tsize\bigcup \{(i-1)_\alpha i_\beta\mid  E_{\alpha\beta}= 2\}.
$$ 
We shall say that a colored partition $\pi$ satisfies the difference $\Cal D$
condition  if $\pi$ does not contain any colored partition from
$\Cal D$, or, equivalently, if $\pi$ is an element of the partition
ideal $\Cal P_\Cal D=\Cal P(A)\backslash (\Cal P(A)\Cal D)$. So, for example,
$\pi=(-5)_1(-3)_8(-2)_9$ does not satisfy the difference $\Cal D$
condition since $E_{89}=2$ and $(-3)_8(-2)_9\in\Cal D$.

Now we can state the following Rogers-Ramanujan type identity:

\proclaim{Theorem 2.1}
$$
\sum_{n=0}^\infty \,\sharp\{|\pi|=n\mid\pi\in\Cal P_\Cal D\}\,q^n=
\prod_{r =1}^\infty \,(1 -q^r)^{-1}.
$$
\endproclaim

The product side follows from the principally specialized Weyl-Kac character
formula for the basic $\frak{sl}(3,\C)\sptilde$-module
\cite{L}, the sum side follows from the
crystal base character formula \cite{KMN${}^2$}. The proof is given in the next
section.

\heading 3. The principally specialized character for the basic
$A_2^{(1)}$-module
\endheading

Let us consider a multiple of the character of
the basic $\frak{sl}(3,\C)\sptilde$-module $L(\Lambda_0)$ 
$$
\prod_{r =1}^\infty \,(1 -e^{-r\delta})^{-1}\cdot e^{-\Lambda_0}\ch L(\Lambda_0).
\tag{3.1}
$$
Then the principal specialization $e^{-\alpha_i}\mapsto q$, $i=0, 1, 2$,
of this product gives
$$
\prod_{r =1}^\infty \,(1 -q^{3r})^{-1}
\prod_{r \not\equiv 0 \mod{3}} \,(1 -q^r)^{-1}
=\prod_{r =1}^\infty \,(1 -q^r)^{-1}.
$$
Here we use notions, notation and results as in \cite{K} or \cite{L}.

On the other side, we may use the results in \cite{KMN${}^2$} to express
(3.1) in terms of colored partitions. We consider a perfect crystal $\Gamma$
for $\frak{sl}(3,\C)\sptilde$ coming from the tensor product of the vector
representation and its dual:
$$\gather
\CD
1  @>1>> 2   @>2>> 5  \\
 @V2VV  @.  @V2VV   \\
3   @>1>> 6   @. 8  \\
 @.  @V1VV  @V1VV  \\
4  @. 7   @>2>> 9 
\endCD \\
\text{together with $0$-arrows}\\
9 @>0>> 4  @>0>> 1 \ ,
\quad 8 @>0>> 3 \quad\text{and}\quad 7 @>0>> 2 \ .
\endgather$$
For $\beta\in\Gamma$ let us denote by $\wt (\beta)$ the $\goth h$-weight of $\beta$,
that is, the restriction of the classical weight of $\beta$ on the fixed Cartan
subalgebra $\goth h\subset\goth g$.
Note that 
$$
\align
&\wt 1=-\wt 9=\alpha_1+\alpha_2,\\
&\wt 2=-\wt 8=\alpha_2,\\
&\wt 3=-\wt 7=\alpha_1,\\
&\wt 4=\wt 5=\wt 6=0.
\endalign
$$
The crystal $\Gamma$ has the unique
energy function $H$ with values in $\{0,1,2\}$, and we have chosen
$$
E_{\alpha\beta}=H(\beta\otimes\alpha).
$$
In particular, we have
$$
E_{44}=H(4\otimes 4)=0,\qquad \wt 4=0.\tag{3.2}
$$
The ground state path for the basic module $L(\Lambda_0)$ is
$$
p_{\Lambda_0}=\, 4,4,4,4,4,4,\dots \ .
$$

By \cite{KMN${}^2$, Proposition 4.6.4}, the set $\Cal P(\Lambda_0,\Gamma)$
of sequences 
$p=\bigl(p(j)\,;\, j\geqslant 1\bigr)$ in $\Gamma$ 
such that $p(j)=p_{\Lambda_0}(j)=4$ for $j\gg 0$ parameterizes a basis of 
$L(\Lambda_0)$, where, by taking into account (3.2), the weight and the degree of a 
sequence is given by
$$\allowdisplaybreaks\align
&|p|= -\sum_{j=1}^\infty j\,H\bigl(p(j+1)\otimes p(j)\bigr), \tag{3.3}\\
&\wt (p) = \sum_{j=1}^\infty \wt \bigl(p(j)\bigr). \tag{3.4}
\endalign$$

We want to interpret this result in terms of colored partitions, and, in
order to do that, let us think of colored partitions in a different way:
Let $\preccurlyeq$ be an order on $\Gamma$ defined as
$$
1\succ 2\succ 3\succ 4\succ 5\succ 6\succ 7\succ 8\succ 9,
$$
and define an order on $A$ by
$$
i_\beta \preccurlyeq j_\gamma \qquad \text{if \qquad  either} \quad  i<j   
\quad \text{or} \quad i=j,\, \beta\preccurlyeq\gamma.
$$ 
This is a total order on $A$.
For a colored partition 
$$
\nu = (j_1)_{\beta_1}\dots (j_s)_{\beta_s}\in \Cal P(A)
$$
 we may assume that  $j_1\leqslant\dots\leqslant j_s< 0$ 
 and that $j_r = j_{r+1}$ implies $\beta_r \preccurlyeq \beta_{r+1}$,
i.e., $(j_r)_{\beta_r} \preccurlyeq (j_{r+1})_{\beta_{r+1}}$.
Sometimes we shall denote a colored partition $\nu$ as 
$$
\nu = \bigl( (j_1)_{\beta_1} \preccurlyeq (j_2)_{\beta_2} \preccurlyeq 
\dots\preccurlyeq (j_s)_{\beta_s}\bigr). 
$$
We may visualize $\nu$ by its Young diagram:
\def\sq{\lower.3ex\vbox{\hrule\hbox{\vrule height1.2ex depth1.2ex\kern2.4ex 
\vrule}\hrule}\,}
$$\spreadlines{-1.15ex}\alignat4
\text{color}\quad &\beta_1\quad &&\sq\sq&&\!\cdots\sq\!\cdots\sq\sq\quad&&-j_1\quad
\text{boxes}\\
&\beta_2\quad&&\sq\sq&&\!\cdots\sq\!\cdots\sq\quad&&-j_2\\
&\quad&&&&\!\cdots &&\\
&\beta_s&&\sq\sq&&\!\cdots\sq &&-j_s \ .
\endalignat
$$
Then the total number of boxes $||\nu||$ in the Young diagram of $\nu$ is
$$
||\nu||=\sum_{m=1}^s \, -j_m\ .
$$
Since to each color $\beta\in\Gamma$ we can associate its weight 
$\wt (\beta)\in\goth h^*$, we define a weight $\wt (\nu)$ of the
colored partition $\nu$ as
$$
\wt (\nu)=\sum_{m=1}^s \, \wt (\beta_m)\ .
$$

Now for a given path  
$$
p=\beta_1,\dots,\beta_{s-1},\beta_s,\dots, 4,4,4,4,\dots
$$ 
we construct a colored partition $\part_\Cal D(p)$ in the
following way: We start with a
large enough $s$, i.e., an $s$ such that $\beta_j=4$ for all $j\geqslant s$,
and we set
$$
-i_s=-i_{s+1}=\dots=0,
$$
and from there on
$$
\align 
&-i_{s-1}= E_{\beta_{s-1}\beta_s},\\
&-i_{s-2}= E_{\beta_{s-2}\beta_{s-1}}+E_{\beta_{s-1}\beta_s}, \\
&\dots\\
&-i_1= \sum^{s-1}_{r=1} E_{\beta_r\beta_{r+1}}.
\endalign
$$
Note that $E_{44}=0$, so it does not matter with which $s\gg 0$ we have started.
Now we define a colored partition $\part_\Cal D(p)$ associated to $p$ as
$$
\part_\Cal D(p)=(-i_1)_{\beta_1}\dots (-i_{s-1})_{\beta_{s-1}}(-i_s)_{\beta_s}\dots
 (0)_4 (0)_4 (0)_4 (0)_4 \dots,
$$
where we identify the product of all $(0)_4$ with $1\in \Cal P(A)$.
It is easy to check that 
$$
E_{\beta_r\beta_{r+1}}=0\quad\text{ implies }\quad
\beta_r\preccurlyeq \beta_{r+1},\tag{3.5}
$$ 
so by construction we have
$$
\part_\Cal D(p)=
\bigl((-i_1)_{\beta_1}\preccurlyeq (-i_2)_{\beta_2}\preccurlyeq\dots\bigr)  .
$$
We may visualize the above construction of $p\mapsto\part_\Cal D(p)$ in
terms of Young diagrams: We start from ``the bottom'' for large enough
$s$ and we add $0$ boxes to a color $\beta_s=4$. As we are finished
with associating $-i_{r+1}$ boxes to $\beta_{r+1}$, we associate to
$\beta_r$ an extra $ E_{\beta_r, \beta_{r+1}}$ boxes, so that
$\part_\Cal D(p)$ satisfies difference conditions
$$
-i_{r}=-i_{r+1}+ E_{\beta_r, \beta_{r+1}}.
$$
By counting the number of boxes we added at each stage, we see that the
total number of boxes in $\part_\Cal D(p)$ equals $-|p|$ given by (3.3), i.e.
$$
||\part_\Cal D(p)||= \sum_{r=1}^s r\,E_{\beta_r, \beta_{r+1}}
=-|p|.\tag{3.6}
$$
We also have (cf. (3.4) and (3.2))
$$
\wt (\part_\Cal D(p))= \sum_{r=1}^s \beta_r=\wt (p).\tag{3.7}
$$

For a given colored partition $\nu$ and a plain partition
$\Delta\in\Cal P(\Z_{<0})$ (i.e. a partition ``without colors''),
$$
\align
&\nu=\bigl((i_1)_{\beta_1}\preccurlyeq\dots \preccurlyeq (i_s)_{\beta_s}\bigr),\\
&\Delta=(j_1\leqslant j_2\leqslant\dots\leqslant j_p),
\endalign
$$
 $i_s, j_p<0$, we define a colored partition $\nu\oplus\Delta$ by
$$
\nu\oplus\Delta=(k_1)_{\beta_1}\dots (k_r)_{\beta_r},
$$
where $r=\max \{s, p\}$, with additional colors 
$\beta_{s+1}=\dots=\beta_p=4$ in the case $p>s$, and
$$
k_n=i_n+\sum_{-j_m\geqslant n}\,-1,
$$
that is, to a Young diagram of $\nu$ we add to each color $\beta_1,\dots,\beta_{-j_1}$
one box, then we add to each color $\beta_1,\dots,\beta_{-j_2}$ another one box,
and so on. Of course, in the case $p>s$ we consider $i_{s+1}=\dots=i_p=0$.
It is clear that in the case $\nu=\part_\Cal D(p)$ we have difference conditions
$$
-i_{r}\geqslant -i_{r+1}+ E_{\beta_r, \beta_{r+1}}.
$$

It is clear that the map 
$$
(p,\Delta)\mapsto \part_\Cal D(p)\oplus \Delta,
\qquad\Cal P(\Lambda_0,\Gamma)\times\Cal P(\Z_{<0})\to\Cal P(A)
$$
 is injective. Hence, by using (3.6) and (3.7), 
\cite{KMN${}^2$, Proposition 4.6.4} implies
$$
\sum_{\pi=\part_\Cal D(p)\oplus \Delta} \,e^{\wt (\pi)-||\pi||\delta}=
\prod_{r =1}^\infty \,(1 -e^{-r\delta})^{-1}\cdot e^{-\Lambda_0}\ch L(\Lambda_0),
$$
where the sum runs over all pats $p\in\Cal P(\Lambda_0,\Gamma)$ 
and plain partitions $\Delta\in\Cal P(\Z_{<0})$.

What we want to see is that every $\pi=\part_\Cal D(p)\oplus \Delta$ is in 
$\Cal P_\Cal D$. For that it is sufficient to show that
$$
(i_\alpha \preccurlyeq j_\beta)\in \Cal P_\Cal D \quad  
\text {if and only if} \quad |i-j|\geqslant E_{\alpha\beta},\tag{3.8}
$$ 
$$
i_\alpha \preccurlyeq j_\beta \preccurlyeq k_\gamma, \  
|i-k| \leqslant 1, \   i_\alpha k_\gamma \in \Cal D \quad\text{implies}\quad  
i_\alpha j_\beta \in \Cal D \quad\text{or} \quad j_\beta k_\gamma \in \Cal D.\tag{3.9}
$$
It is easy to see, by using (3.5), that (3.8) holds. 
Now, in the presence of (3.8), we can easily check 
that (3.9) holds as well; in terms of $E$ it reads as
$$
E_{\alpha\gamma} \leqslant E_{\alpha\beta}+E_{\beta\gamma},\tag{3.10}
$$
or, equivalently, as
$$
E_{\alpha\gamma} \leqslant E_{\alpha\beta}\quad\text{if}\quad E_{\beta\gamma}=0,
\qquad
E_{\alpha\gamma} \leqslant E_{\beta\gamma}\quad\text{if}\quad  E_{\alpha\beta}=0.
$$
Hence every $\pi=\part_\Cal D(p)\oplus \Delta$ is in 
$\Cal P_\Cal D$. Moreover, every $\pi\in\Cal P_\Cal D$ can be written in this
way, and hence we have
$$
\sum_{\pi\in \Cal P_\Cal D} \,e^{\wt (\pi)-||\pi||\delta}=
\prod_{r =1}^\infty \,(1 -e^{-r\delta})^{-1}\cdot e^{-\Lambda_0}\ch L(\Lambda_0).
$$
Since the map $(-i)_\alpha\mapsto |(-i)_\alpha|$ defined by (2.1) is the principal 
specialization, Theorem 2.1 holds.


\heading 4. Some remarks
\endheading

Let $\Gamma$ be a classical crystal with an energy function $H$ with values
in the set $\{0,1,2\}$. Then, as before, we can consider
colored partitions $\pi\in\Cal P(A)$ with $A=\Gamma\times\Z_{<0}$
and we can define the degree $||\pi||$ and the $\goth h$-weight of $\pi$ as above,
with $\goth h=\lspan \{h_1,\dots,h_\ell\}$ (see \cite{KMN${}^2$}). 
If we set $E_{\alpha\beta}=H(\beta\otimes\alpha)\in\{0,1,2\}$, we can define
$$
\Cal D =\{i_\alpha i_\beta\mid E_{\alpha\beta} E_{\beta\alpha}\geqslant 1\}
\tsize\bigcup \{(i-1)_\alpha i_\beta\mid  E_{\alpha\beta}= 2\},
$$
and we can consider colored partitions which satisfy difference $\Cal D$
conditions, that is, $\pi\in\Cal P_\Cal D$. Let us define a
``character'' of the partition ideal $\Cal P_\Cal D$ as
$$
\ch (\Cal P_\Cal D)=\sum_{\pi\in \Cal P_\Cal D} \,e^{\wt (\pi)-||\pi||\delta}.
$$
The map
 $e^{-\alpha_i}\mapsto q$, $i=0, 1,\dots,\ell$, defines $e^{-\delta}\mapsto q^m$
and the principally specialized character
$$
\ch_q (\Cal P_\Cal D)=
\sum_{n=0}^\infty \,\sharp\{|\pi|=n\mid\pi\in\Cal P_\Cal D\}\,q^n,
$$
where $|(-j)_{\beta}|=mj-\la\rho,\beta\ra$, with $\la\rho,\alpha_i\ra=1$
for $i=1,\dots,\ell$ (cf. \cite{K}).

With this notation at hand we can write Theorem 2.1, for our particular
choice of the $A_2^{(1)}$-crystal $\Gamma$, as
$$
\ch_q (\Cal P_\Cal D)=
\prod_{r =1}^\infty \,(1 -q^r)^{-1}.
$$

 Let $\Gamma$ be a perfect crystal
for $\tilde\goth g=\frak{sl}(n,\C)\sptilde$, $n\geq 4$, 
coming from the tensor product of the vector
representation and its dual. Then there is an energy function $H$ taking values in 
$\{0,1,2\}$ and the ground state path for the basic $\tilde\goth g$-module
is a constant sequence $p(j)=a$, with $E_{aa}=H(a\otimes a) =0$ and $\wt (a)=0$.
Moreover, there is a total order $\preccurlyeq$ on $\Gamma$ such that (3.5)
holds. So the same proof would go through, and the above identity would hold, 
if we could show (3.10).
For example, in the case of $A_3^{(1)}$-crystal
$$
\gather
\CD
14  @>1>> 24   @>2>> 34 @>3>> 44  \\
  @V3VV  @V3VV  @.  @V3VV    \\
13  @>1>> 23   @>2>> 33 @. 43  \\
  @V2VV  @.  @V2VV  @V2VV    \\
12  @>1>> 22   @. 32 @>3>> 42  \\
  @.  @V1VV  @V1VV  @V1VV    \\
11  @. 21   @>2>> 31 @>3>> 41 
\endCD \\
\text{together with $0$-arrows}\\
41 @>0>> 11  @>0>> 14 \ , 
\quad 21 @>0>> 24\ , \quad 31 @>0>> 34 \ ,
\quad 42 @>0>> 12\ , \quad 43 @>0>> 13 
\endgather$$
all these properties hold, including (3.10), and we have an identity of
the above form.
On the other side, let $\Gamma=\{1,2,3,4\}$ be the $A_1^{(1)}$-crystal
$$
\CD
1  @>1>> 3  \\
 @AA0A    @V1VV   \\
2   @<<0< 4    
\endCD
$$
with the energy matrix
$$
E =\left(\matrix 
2 & 1& 2& 2\\
1 & 0& 1& 1\\
0 & 1& 0& 2\\
0 & 1& 0& 2
\endmatrix\right) .
$$
The KMN${}^2$ crystal base character formula
is proved under the assumption that the rank of $\goth g$ is at
least two, but still many results also hold for $\tilde\goth g=\frak{sl}(2,\C)\sptilde$.
So it is reasonable to ask whether an analogue of Theorem 2.1 
holds as well, i.e., whether
$$
\ch_q (\Cal P_\Cal D)=
\prod_{r =1}^\infty \,(1 -q^r)^{-1}\quad ?
$$
Note that here the principal specialization reads
$$
(-i)_1 \mapsto {2i -1}, \quad
(-i)_2 \mapsto 2i,  \quad
(-i)_3 \mapsto 2i,  \quad
(-i)_4 \mapsto 2i+1.  
$$

What is surprising is that there are other identities of a similar form, but
which are not related to the crystal base theory, at least not in any obvious
way: consider an  ``almost perfect'' $A_1^{(1)}$-crystal
$$
1 \overset{\overset 1\to\longrightarrow}\to{\underset 0\to\longleftarrow}
2 \overset{\overset 1\to\longrightarrow}\to{\underset 0\to\longleftarrow} 3 
$$
with the energy matrix
$$
E =\left(\matrix 
2 & 2& 2\\
1 & 1& 2\\
0 & 1& 2
\endmatrix\right) \ .
$$
Then we have
$$
\ch_q (\Cal P_\Cal D)=
\prod_{r \text{ odd}} \,(1 -q^r)^{-1}.\tag{4.1}
$$
Note that here the principal specialization reads
$$
(-i)_1 \mapsto {2i -1}, \quad
(-i)_2 \mapsto 2i,  \quad
(-i)_3 \mapsto 2i+1.  
$$
Moreover, if we define a map $(-i)_\alpha\mapsto |(-i)_\alpha|$ for $i>0$ by
$$
(-i)_1 \mapsto {3i -2}, \quad
(-i)_2 \mapsto 3i,  \quad
(-i)_3 \mapsto 3i+2,  
$$
(i.e. if we take the (1,2)-specialization of $\ch (\Cal P_\Cal D)$),
then we have a Capparelli identity (see \cite{C1}--\cite{C3}, \cite{A2})
$$
\sum_{n=0}^\infty \,\sharp\{|\pi|=n\mid\pi\in\Cal P_\Cal D\}\,q^n=
 \prod\Sb r \equiv 1,3,5,6 \mod{6}\endSb
 \,(1 + q^r).\tag{4.2}
$$
Both (4.1) and (4.2) are proved in \cite{MP2} as specializations of 
the character formula for the basic $\frak{sl}(2,\C)\sptilde$-module
written in the form
$$
e^{-\Lambda_0} \ch\, L(\Lambda_0) = \ch (\Cal P_\Cal D).
$$
The character formula itself is proved by using the Lepowsky-Wilson
approach, the proof being quite parallel to \cite{MP1}. 
The set $\Cal D$ is originally defined as the set of
leading terms for the
vertex operator algebra defining relations for the basic 
$\frak{sl}(2,\C)\sptilde$-module (cf. \cite{MP2, Section 6}). 

So it seems that interesting combinatorial properties of $\ch (\Cal P_\Cal D)$
go beyond the crystal base character formula \cite{KMN${}^2$, Proposition 4.6.4},
at least when the relation
$E_{\alpha\gamma} \leqslant E_{\alpha\beta}+E_{\beta\gamma}$ holds
and when there is a total order $\preccurlyeq$ on $\Gamma$ such that
$E_{\alpha\beta}=0$ implies $\alpha\preccurlyeq\beta$. Further
indications for this provide the results in \cite{P2}
and the examples below.

Although the one by one matrices have nothing to do with crystals, the
notion of difference $\Cal D$ condition still makes sense for 
$\Gamma=\{1\}$ and $E=(E_{11})=(2)$;
it is simply the difference $2$ condition and
$\ch_q (\Cal P_\Cal D)$, defined via $|(-i)_1|=i$,
 is the sum side of a Rogers-Ramanujan identity.
 On the other hand, for $E=(E_{11})=(1)$
the difference $\Cal D$ condition defines partitions in distinct parts.
Of course, the second case is much simpler, and it has an equally simple
analogue for $\Gamma=\{1, 2\}$ and 
$
E =\left(\smallmatrix 
1 & 1\\
0 & 1
\endsmallmatrix\right)  
$:
consider an  ``almost perfect'' $A_1^{(1)}$-crystal
$$
1 \overset{\overset 1\to\longrightarrow}\to{\underset 0\to\longleftarrow} 2 \ 
$$
with the energy matrix $E_{\alpha\beta}=H(\beta\otimes\alpha)$ chosen to be $E$. \,
Then 
$$
(-j)_1\mapsto 2j-\tfrac12,\qquad (-j)_2\mapsto 2j+\tfrac12
$$ 
is the principal
specialization and we have an identity for partitions in half-integers
$$
\ch_q (\Cal P_\Cal D)=
\prod_{r \geqslant 1} \,(1 +q^{r+\frac12}).
$$
In the case 
$
E =\left(\smallmatrix 
2 & 2\\
1 & 2
\endsmallmatrix\right)  
$
the principal specialization of $\Cal P_\Cal D$
gives partitions in half-integers satisfying 
difference 3 condition, but I am not aware of any formula that would express
$\ch_q (\Cal P_\Cal D)$ as an infinite product.

Finally, the formulation of Rogers-Ramanujan type identities for colored
partitions in terms of energy functions for crystals was also
motivated by a desire to understand better an identity for
the basic $A_2^{(1)}$-module obtained in \cite{MP3}.
The set $\Cal D$ of difference conditions is defined as 
the set of leading terms of relations
for the basic module, but can be defined as well in the following way:
consider the weighted $A_2$-crystal \ $\Gamma '=\Gamma\backslash \{4\}$
$$
\CD
1  @>1>> 2   @>2>> 5  \\
 @V2VV  @.  @V2VV   \\
3   @>1>> 6   @. 8  \\
 @.  @V1VV  @V1VV  \\
  @. 7   @>2>> 9 
\endCD \quad . 
$$
It is not possible do define $0$-arrows which would turn $\Gamma'$ into
an $A_2^{(1)}$-crystal with an energy function. But still, we can
define a ``difference condition matrix'' $(E_{\alpha\beta})_{\alpha,\beta\in\Gamma'}$
$$
E =\left(\matrix 
2 & 2& 2&  2& 2& 2& 2& 2\\
1 & 2& 1&  2& 1& 2& 2& 2\\
1 & 1& 2&  1& 2& 2& 2& 2\\
0 & 1& 1&  1& 1& 1& 2& 2\\
1 & 1& 1&  1& 1& 2& 1& 2\\
0 & 1& 0&  1& 1& 2& 1& 2\\
0 & 0& 1&  1& 1& 1& 2& 2\\
0 & 0& 0&  0& 1& 1& 1& 2 
\endmatrix\right) 
$$
and we can define a subset $\Cal D$ in $\Cal P(A)$, $A=\Gamma'\times\Z_{<0}$, by
$$\align
\Cal D =&\{i_\alpha i_\beta\mid E_{\alpha\beta} E_{\beta\alpha}\geqslant 1\}
\tsize\bigcup \{(i-1)_\alpha i_\beta\mid  E_{\alpha\beta}= 2\}\\
&\tsize\bigcup\{(i-1)_3\, i_5\, i_1\}\tsize\bigcup\{ (i-1)_9\,(i-1)_5\, i_7\}.
\endalign$$
Then, as proved in \cite{MP3}, we have
$$
e^{-\Lambda_0} \ch\, L(\Lambda_0) = \ch (\Cal P_\Cal D).
$$
By taking the principal specialization  (2.1) without the color
$4$, we obtain a combinatorial identity
$$
\ch_q (\Cal P_\Cal D)=
\prod_{r \not\equiv 0 \mod{3}} \,(1 -q^r)^{-1}.
$$
 As it happens, the above ``difference 2 conditions'' 
$\{(i-1)_\alpha i_\beta\mid  E_{\alpha\beta}= 2\}$ are the same
as in the case of partitions discussed in Section 2. One is 
tempted to think that this is more than a mere coincidence and that
the difference conditions
$\{(i-1)_3\, i_5\, i_1\}\tsize\bigcup\{ (i-1)_9\,(i-1)_5\, i_7\}$
are some sort of ``corrections'' of the fact that $\Gamma'$ is
not a perfect $A_2^{(1)}$-crystal. 
With this regard it may be interesting to note that the energy matrix
$E$ in Section 2 is invariant under the change
$$
1\leftrightarrow 1, \quad 2\leftrightarrow 3, \quad 5\leftrightarrow 6, \quad 
7\leftrightarrow 8, \quad 9\leftrightarrow 9
$$
and $4\leftrightarrow 4$. Here we have 
$E_{51}=0$, $E_{61}=1$ and $E_{95}=0$, $E_{96}=1$. So the
``interaction'' between parts $i_6\, i_1\in \Cal D$ and
$i_5\, i_1\not\in \Cal D$ is not symmetrical, the later is ``compensated''
with a weaker requirement $(i-1)_3\, i_5\, i_1\in \Cal D$. Likewise
$i_9\, i_6\in \Cal D$ ``corresponds'' to a weaker $i_9\, i_5\, (i+1)_7\in \Cal D$.
It should be said that the present proof
of the above identity has no 
connections with the crystal base theory.


\enddocument